\documentclass[10pt, a4paper]{article}
\usepackage{latexsym}
\usepackage{amsmath}
\usepackage{amssymb}
\usepackage[all]{xy}
\usepackage{bm}
\usepackage{amsfonts}
\usepackage{verbatim}
\usepackage{amsthm}
\usepackage{mathrsfs}
\usepackage{epsfig}
\usepackage{xy}
\usepackage{array}
\usepackage{stmaryrd}
\usepackage{graphicx,color}
\usepackage{xcolor}
\usepackage{tikz}
\usetikzlibrary{arrows,calc}
\usepackage{etex}
\usepackage{mathdots}
\usepackage{float}
\usepackage{graphics}
\usepackage{pdflscape}
\usepackage{cite}

\usepackage{anysize,hyperref}
\input xypic
\xyoption{all}

\usepackage[perpage,symbol]{footmisc}
\usepackage{setspace}
\renewcommand{\cal}{\mathcal}

\def\C{\mathscr{C}}
\def\E{\mathbb{E}}
\def\s{\mathfrak{s}}

\def\op{^\mathrm{op}}
\def\Ab{\mathit{Ab}}
\def\del{\delta}
\def\dr{\ar@{->}[r]}
\def\X{\mathscr{X}}

\def\X{\mathscr{X}}

\def\Hom{\mbox{Hom}}

\begin{document}

\baselineskip=15pt
\title{\Large{\bf Grothendieck groups in extriangulated categories}}
\medskip
\author{\textbf{Bin Zhu and Xiao Zhuang}}

\date{}

\maketitle
\def\blue{\color{blue}}
\def\red{\color{red}}

\newtheorem{theorem}{Theorem}[section]
\newtheorem{lemma}[theorem]{Lemma}
\newtheorem{corollary}[theorem]{Corollary}
\newtheorem{proposition}[theorem]{Proposition}
\newtheorem{conjecture}{Conjecture}
\theoremstyle{definition}
\newtheorem{definition}[theorem]{Definition}
\newtheorem{question}[theorem]{Question}
\newtheorem{remark}[theorem]{Remark}
\newtheorem{notation}[theorem]{Notation}
\newtheorem{remark*}[]{Remark}
\newtheorem{example}[theorem]{Example}
\newtheorem{example*}[]{Example}

\newtheorem{construction}[theorem]{Construction}
\newtheorem{construction*}[]{Construction}

\newtheorem{assumption}[theorem]{Assumption}
\newtheorem{assumption*}[]{Assumption}

\baselineskip=17pt
\parindent=0.5cm

\begin{abstract}

The aim of the paper is to discuss the relation subgroups of the Grothendieck groups of extriangulated categories and certain other subgroups. It is shown that a locally finite
extriangulated category $\C$ has Auslander-Reiten $\E-$triangles and the relations of the Grothendieck group $K_{0}(\C)$ are
generated by the Auslander-Rieten $\E-$triangles. A partial converse result is given when restricting to the triangulated categories with a cluster tilting subcategory: in the triangulated category $\C$ with a cluster tilting subcategory, the relations of the Grothendieck group $K_0(\C)$ are generated by Auslander-Reiten triangles if and only if the triangulated category $\C$ is locally finite. It is also shown that there is a one-to-one correspondence between subgroups of $K_{0}(\C)$ containing the image of $\mathcal G$ and dense $\mathcal G-$(co)resolving subcategories of $\C$ where $\mathcal G$ is a generator of $\C,$ which generalizes results about classifying subcategories of a triangulated \cite{t} or an exact category $\C$ \cite{m} by subgroups of $K_{0}(\C)$. \\[0.5cm]

\medskip
\end{abstract}

\section{Introduction}
 Auslander-Reiten sequences were introduced by Auslander and Reiten in \cite{ar1} as a useful tool to study the local structure of the
 module category of an artin algebra. Since then, it has been generalized to the situation of exact categories [\cite{asm}, \cite{k}],
 triangulated categories \cite{h} and its subcategories [\cite{asm}, \cite{j}] and some additive categories [\cite{liu}, \cite{sh}, \cite{i}] by many authors. Recently, extriangulated categories were introduced by Nakaoka and Palu \cite{np}
 to give a simultaneous generalization of exact categories and triangulated categories. As in the situation of exact categories or
 triangulated categories, there is also a notion of Auslander-Reiten $\E-$triangles in extriangulated categories, which was first introduced
 in \cite{zz}, and then was fully developed in \cite{inp}. For an extriangulated category $\C$, one can define its Grothendieck group $K_0(\C)$ similarly as the case of exact categories or triangulated categories: for an object $X$ in $\C$, we use $[X]$ to denote the isoclass of $X$, and $F$ to denote the free group generated by the isoclasses $[X]$ of objects $X$ in $\C$, then the Grothendieck group $K_0(\C)$ of $\C$ is defined as the quotient of $F$ by the subgroup generated by $[\delta]=[A]+[C]-[B]$ for all $\E-$triangles  $\xymatrix{A\ar[r]^{}&B\ar[r]^{}&C\ar@{-->}[r]^{\delta}&}$ in $\C,$ (For more detail, see Section 4, or Part II in \cite{pppp}]. All such $[\delta]$, or the subgroup generated by them are called the relations or the relation group of $K_0(\C)$ respectively. If $\C$ is the module category mod$\Lambda$ of an artin algebra $\Lambda$, there is a well-known result [\cite{b}, \cite{a}] saying that the relations of $K_{0}(\C)$ are generated by Auslander-Reiten sequences if and only if $\Lambda$ is of finite representation type. Similarly, when $\C$ is a locally finite triangulated
 category, Xiao and Zhu [\cite{xz}, Theorem 2.1] showed that the Auslander-Reiten triangles generate the relations for the Grothendieck group.
 Beligiannis [\cite{be}, Theorem 12.1] proved the converse of this result holds when $\C$ is a compactly generated triangulated category. Recently,
 Haugland [\cite{ha1}, Theorem 2.3] and [\cite{pppp}, Theorem 3.8, Theorem 4.13] showed the converse of the result in some special cases.
 In this paper, we show that locally finite extriangulated categories have Auslander-Reiten $\E-$triangles and that for a locally finite extriangulated category, the relations of the Grothendieck group are generated by
 Auslander-Reiten $\E-$triangles. The later result is a generalization of results in module categories [\cite{a}, \cite{ar2}, \cite{b}] and triangulated
 categories \cite{xz}. We give a partial converse result when restricting to the triangulated category with a cluster tilting subcategory. We show that in these triangulated categories, the relations of the Grothendieck group are generated by Auslander-Reiten triangles if and only if these triangulated categories are locally finite. At the end of the paper, we consider the relations between subgroups of the Grothendieck group $K_{0}(\C)$ and
 subcategories of the extriangulated category $\C.$ We show that for an extriangulated category with a generator $\mathcal G,$ there is a
 one-to-one correspondence between subgroups of $K_{0}(\C)$ containing the image of $\mathcal G$ and dense $\mathcal G-$(co)resolving
 subcategories of $\C.$ This can be seen as a generalization of the results in triangulated categories \cite{t} and exact
 categories \cite{m}. Independently of our work, J. Haugland has recently generalized the above result to the situation of $n-$exangulated categories in \cite{ha2}.

 The paper is organized as follows. In Section 2, we recall the definition of an extriangulated category and outline some basic properties
 that will be used later. In Section 3 we recall the definition of Auslander-Reiten $\E-$triangles in an extriangulated category \cite{inp} and show that every locally
 finite extriangulated category has Auslander-Reiten $\E-$triangles. In Section 4 we show that
 for a locally finite extriangulated category, the relations of the Grothendieck group are all generated by Auslander-Reiten $\E-$triangles and we give a partial converse result for the triangulated categories with a cluster tilting subcategory. In section 5 we show that for an extriangulated category $\C$ with a generator $\mathcal G,$ there is a one-to-one correspondence between subgroups of $K_{0}(\C)$ containing the image of $\mathcal G$ and dense $\mathcal G-$(co)resolving subcategories of $\C.$

\section{Preliminaries}
Throughout this article, $\C$ denotes a K-linear Krull-Schmidt category with K is a field, and all subcategories are full additive subcategories
closed under isomorphisms.
We denote by $\C(A,B)$ or $\Hom_{\C}(A,B)$ the set of morphisms from $A$ to $ B$ in the category $\C$.

Recall the definition of functorially finite subcategories in $\C$.
Let $\X$ be a subcategory of $\C$, and $B\in \C$. A
morphism $f_B\colon X_B\to B$ of $\C$ with $X_B$ an object in $\X$ ($g_B\colon B\to Y_B$ of $\C$ with $Y_B$ an object in $\X$, resp.)
is said to be a {right $\X$-approximation} (left $\X$-approximation, resp.) of $B$, if the
morphisms $\C(X,f_B)\colon\C(X,X_B)\to\C(X,B)$ ($\C(g_B,X)\colon\C(B,X)\to\C(Y_B,X)$, resp.) are surjectives
for all objects $X$ in $\X$.
Moreover, one says that
$\X$ is contravariantly finite (covariantly finite, resp.) in $\C$ if every object in $\C$ has a right (left, resp.) $\X$-approximation. $\X$ is functorially finite if it is both contravariantly finite and covariantly finite in $\C$.

Now we outline definitions and some basic properties of extriangulated categories from \cite{np}.

Let $\C$ be an additive category. Suppose that $\C$ is equipped with a biadditive functor $\E\colon\C\op\times\C\to\Ab$. For any pair of
objects $A,C\in\C$, an element $\delta\in\E(C,A)$ is called an $\E$-extension. Thus formally, an $\E$-extension is a triplet $(A,\delta,C)$.
Let $(A,\del,C)$ be an $\E$-extension. Since $\E$ is a bifunctor, for any $a\in\C(A,A')$ and $c\in\C(C',C)$, we have $\E$-extensions
$$ \E(C,a)(\del)\in\E(C,A')\ \ \text{and}\ \ \ \E(c,A)(\del)\in\E(C',A). $$
We abbreviate denote them by $a_\ast\del$ and $c^\ast\del$ respectively.
For any $A,C\in\C$, the zero element $0\in\E(C,A)$ is called the spilt $\E$-extension.

\begin{definition}{[\cite{np}, Definition 2.3]}\label{mo}
Let $(A,\del,C),(A',\del',C')$ be any pair of $\E$-extensions. A morphism $$(a,c)\colon(A,\del,C)\to(A',\del',C')$$ of $\E$-extensions is a
pair of morphisms $a\in\C(A,A')$ and $c\in\C(C,C')$ in $\C$, satisfying the equality
$$ a_\ast\del=c^\ast\del'. $$
Simply we denote it as $(a,c)\colon\del\to\del'$.
\end{definition}

Let $A,C\in\C$ be any pair of objects. Sequences of morphisms in $\C$
$$\xymatrix@C=0.7cm{A\ar[r]^{x} & B \ar[r]^{y} & C}\ \ \text{and}\ \ \ \xymatrix@C=0.7cm{A\ar[r]^{x'} & B' \ar[r]^{y'} & C}$$
are said to be equivalent if there exists an isomorphism $b\in\C(B,B')$ which makes the following diagram commutative.
$$\xymatrix{
A \ar[r]^x \ar@{=}[d] & B\ar[r]^y \ar[d]_{\simeq}^{b} & C\ar@{=}[d]&\\
A\ar[r]^{x'} & B' \ar[r]^{y'} & C &}$$

We denote the equivalence class of $\xymatrix@C=0.7cm{A\ar[r]^{x} & B \ar[r]^{y} & C}$ by $[\xymatrix@C=0.7cm{A\ar[r]^{x} & B \ar[r]^{y} &
C}]$.
For any $A,C\in\C$, we denote as
$ 0=[A\xrightarrow{\binom{1}{0}}A\oplus C\xrightarrow{(0,\ 1)}C].$

\begin{definition}{[\cite{np}, Definition 2.9]}\label{re}
For any $\E$-extension $\delta\in\E(C,A)$, one can associate an equivalence class $\s(\delta)=[\xymatrix@C=0.7cm{A\ar[r]^{x} & B \ar[r]^{y} &
C}].$  This $\s$ is called a realization of $\E$, if it satisfies the following condition:
\begin{itemize}
\item Let $\del\in\E(C,A)$ and $\del'\in\E(C',A')$ be any pair of $\E$-extensions, with $$\s(\del)=[\xymatrix@C=0.7cm{A\ar[r]^{x} & B
    \ar[r]^{y} & C}],\ \ \ \s(\del')=[\xymatrix@C=0.7cm{A'\ar[r]^{x'} & B'\ar[r]^{y'} & C'}].$$
Then, for any morphism $(a,c)\colon\del\to\del'$, there exists $b\in\C(B,B')$ which makes the following diagram commutative.
\begin{equation}\label{t1}
\begin{array}{l}
$$\xymatrix{
A \ar[r]^x \ar[d]^a & B\ar[r]^y \ar[d]^{b} & C\ar[d]^c&\\
A'\ar[r]^{x'} & B' \ar[r]^{y'} & C' &}$$
\end{array}
\end{equation}
\end{itemize}
In this case, we say that sequence $\xymatrix@C=0.7cm{A\ar[r]^{x} & B \ar[r]^{y} & C}$ realizes\ $\del$, whenever it satisfies
$\s(\del)=[\xymatrix@C=0.7cm{A\ar[r]^{x} & B \ar[r]^{y} & C}]$.
Remark that this condition does not depend on the choices of the representatives of the equivalence classes.
In the above situation, we say that (\ref{t1}) (or the triplet $(a,b,c)$) realizes $(a,c)$.
\end{definition}

\begin{definition}{[\cite{np}, Definition 2.10]}\label{are}
Let $\C, \E$ be as above. A realization $\s$ of $\E$ is said to be additive if the following conditions are satisfied:
\begin{itemize}
\item[{\rm (i)}] For any $A,C\in\C$, the split $\E-$extension $0\in \E(C,A)$ satisfies $\s(0)=0$.
\item[{\rm (ii)}] For any pair of $\E-$extensions $\del=(A,\del,C)$ and $\del'=(A',\del',C')$,
                  $$\s(\del\oplus \del')=\s(\del)\oplus \s(\del')$$ holds.
\end{itemize}
\end{definition}

\begin{definition}{[\cite{np}, Definition 2.12]}\label{ex}
We call the pair $(\E,\s)$ an external triangulation of $\C$ if the following conditions are satisfied:
\begin{itemize}
\item[{\rm (ET1)}] $\E\colon\C\op\times\C\to\Ab$ is a biadditive functor.
\item[{\rm (ET2)}] $\s$ is an additive realization of $\E$.
\item[{\rm (ET3)}] Let $\del\in\E(C,A)$ and $\del'\in\E(C',A')$ be any pair of $\E$-extensions, realized as
$$ \s(\del)=[\xymatrix@C=0.7cm{A\ar[r]^{x} & B \ar[r]^{y} & C}],\ \ \s(\del')=[\xymatrix@C=0.7cm{A'\ar[r]^{x'} & B' \ar[r]^{y'} & C'}]. $$
For any commutative square
$$\xymatrix{
A \ar[r]^x \ar[d]^a & B\ar[r]^y \ar[d]^{b} & C&\\
A'\ar[r]^{x'} & B' \ar[r]^{y'} & C' &}$$
in $\C$, there exists a morphism $(a,c)\colon\del\to\del'$ which is realized by $(a,b,c)$.
\item[{\rm (ET3)$\op$}] Let $\del\in\E(C,A)$ and $\del'\in\E(C',A')$ be any pair of $\E$-extensions, realized by
$$\xymatrix@C=0.7cm{A\ar[r]^{x} & B \ar[r]^{y} & C}\ \ \text{and}\ \ \ \xymatrix@C=0.7cm{A'\ar[r]^{x'} & B' \ar[r]^{y'} & C'}$$
respectively.
For any commutative square
$$\xymatrix{
A \ar[r]^x& B\ar[r]^y \ar[d]^{b} & C\ar[d]^c&\\
A'\ar[r]^{x'} & B' \ar[r]^{y'} & C' &}$$
in $\C$, there exists a morphism $(a,c)\colon\del\to\del'$ which is realized by $(a,b,c)$.

\item[{\rm (ET4)}] Let $(A,\del,D)$ and $(B,\del',F)$ be $\E$-extensions realized by
$$\xymatrix@C=0.7cm{A\ar[r]^{f} & B \ar[r]^{f'} & D}\ \ \text{and}\ \ \ \xymatrix@C=0.7cm{B\ar[r]^{g} & C \ar[r]^{g'} & F}$$
respectively. Then there exist an object $E\in\C$, a commutative diagram
$$\xymatrix{A\ar[r]^{f}\ar@{=}[d]&B\ar[r]^{f'}\ar[d]^{g}&D\ar[d]^{d}\\
A\ar[r]^{h}&C\ar[d]^{g'}\ar[r]^{h'}&E\ar[d]^{e}\\
&F\ar@{=}[r]&F}$$
in $\C$, and an $\E$-extension $\del^{''}\in\E(E,A)$ realized by $\xymatrix@C=0.7cm{A\ar[r]^{h} & C \ar[r]^{h'} & E},$ which satisfy the
following compatibilities.
\begin{itemize}
\item[{\rm (i)}] $\xymatrix@C=0.7cm{D\ar[r]^{d} & E \ar[r]^{e} & F}$  realizes $f'_{\ast}\del'$,
\item[{\rm (ii)}] $d^\ast\del''=\del$,

\item[{\rm (iii)}] $f_{\ast}\del''=e^{\ast}\del'$.
\end{itemize}

\item[{\rm (ET4)$\op$}]  Let $(D,\del,B)$ and $(F,\del',C)$ be $\E$-extensions realized by
$$\xymatrix@C=0.7cm{D\ar[r]^{f'} & A \ar[r]^{f} & B}\ \ \text{and}\ \ \ \xymatrix@C=0.7cm{F\ar[r]^{g'} & B \ar[r]^{g} & C}$$
respectively. Then there exist an object $E\in\C$, a commutative diagram
$$\xymatrix{D\ar[r]^{d}\ar@{=}[d]&E\ar[r]^{e}\ar[d]^{h'}&F\ar[d]^{g'}\\
D\ar[r]^{f'}&A\ar[d]^{h}\ar[r]^{f}&B\ar[d]^{g}\\
&C\ar@{=}[r]&C}$$
in $\C$, and an $\E$-extension $\del^{''}\in\E(C,E)$ realized by $\xymatrix@C=0.7cm{E\ar[r]^{h'} & A \ar[r]^{h} & C},$ which satisfy the
following compatibilities.
\begin{itemize}
\item[{\rm (i)}] $\xymatrix@C=0.7cm{D\ar[r]^{d} & E \ar[r]^{e} & F}$  realizes $g'^{\ast}\del$,
\item[{\rm (ii)}] $\del'=e_\ast\del''$,

\item[{\rm (iii)}] $d_\ast\del=g^{\ast}\del''$.
\end{itemize}
\end{itemize}
In this case, we call $\s$ an $\E$-triangulation of $\C$, and call the triplet $(\C,\E,\s)$ an externally triangulated category, or for short,
extriangulated category.
\end{definition}

For an extriangulated category $\C$, we use the following notation:

\begin{itemize}
\item A sequence $\xymatrix@C=0.41cm{A\ar[r]^{x} & B \ar[r]^{y} & C}$ is called a conflation if it realizes some $\E$-extension
    $\del\in\E(C,A)$.
\item A morphism $f\in\C(A,B)$ is called an inflation if it admits some conflation $\xymatrix@C=0.7cm{A\ar[r]^{f} & B \ar[r]& C}.$
\item A morphism $f\in\C(A,B)$ is called a deflation if it admits some conflation $\xymatrix@C=0.7cm{K\ar[r]& A \ar[r]^f& B}.$
\item If a conflation $\xymatrix@C=0.6cm{A\ar[r]^{x} & B \ar[r]^{y} & C}$ realizes $\del\in\E(C,A)$, we call the pair
    $(\xymatrix@C=0.41cm{A\ar[r]^{x} & B \ar[r]^{y} & C},\del)$ an $\E$-triangle, and write it in the following way.
$$\xymatrix{A\ar[r]^{x} & B \ar[r]^{y} & C\ar@{-->}[r]^{\del}&}$$
\item Given an $\E-$triangle $\xymatrix{A\ar[r]^{x} & B \ar[r]^{y} & C\ar@{-->}[r]^{\del}&},$ we call $A$ the CoCone of $y\colon B\to C,$
    and denote it by CoCone$(B\to C),$ also denote by CoCone$(y);$ we call $C$ the Cone of $x\colon A\to B,$ and denote it by Cone$(A\to
    B),$ also denote by Cone$(x).$
\item Let $\xymatrix{A\ar[r]^{x}&B\ar[r]^{y}&C\ar@{-->}[r]^{\delta}&}$ and
    $\xymatrix{A'\ar[r]^{x'}&B'\ar[r]^{y'}&C'\ar@{-->}[r]^{\delta'}&}$ be any pair of $\E$-triangles. If a triplet $(a,b,c)$ realizes
    $(a,c)\colon\del\to\del'$ as in $(\ref{t1})$, then we write it as
    $$\xymatrix{
A \ar[r]^x \ar[d]_{a} & B\ar[r]^y \ar[d]^{b} & C\ar@{-->}[r]^{\delta} \ar[d]^{c}&\\
A'\ar[r]^{x'} & B' \ar[r]^{y'} & C' \ar@{-->}[r]^{\delta'}&}$$
and call $(a,b,c)$ a morphism of $\E$-triangles.
\item A subcategory $\cal{T}$ of $\C$ is called extension-closed if $\cal{T}$ is closed under extensions, i.e. for any $\E$-triangle
    $\xymatrix{A\ar[r]^{x} & B \ar[r]^{y} & C\ar@{-->}[r]^{\del}&} $ with $A, C\in \cal{T}$, we have $B\in \cal{T}$.
\end{itemize}

\begin{remark}\label{ze}
For any pair of objects $A,B\in\C$, we have  $\E$-triangles in $\C$: $$\xymatrix{A\ar[r]^{\binom{1}{0}\quad}&A\oplus B\ar[r]^{\quad(0,\
1)}&B\ar@{-->}[r]^{0}&}\ \ and \ \ \xymatrix{A\ar[r]^{\binom{0}{1}\quad}&B\oplus A\ar[r]^{\quad(1,\, 0)}&B\ar@{-->}[r]^{0}&}.$$
\end{remark}

\begin{example}\label{ext}
(1) Exact categories and extension-closed subcategories of triangulated categories are extriangulated categories. The extension-closed
subcategories of an extriangulated category are extriangulated categories, see \cite{np} for more detail.
\vspace{2mm}

(2) Let $\C$ be an extriangulated category, and $\cal J$ be a subcategory of $\C$.
If $\cal J\subseteq\cal P\cap\cal I$, where $\cal P$ is the subcategory of projective objects in $\C$ and $\cal I$ is the subcategory of
injective objects in $\C$ (see Definition 2.10), then $\C/\cal J$ is an extriangulated category.
This construction gives extriangulated categories which are neither exact nor triangulated in general.
For more details, see [\cite{np}, Proposition 3.30]. There are also other such examples in \cite{zz}.
\end{example}

There are some basic results on extriangulated categories which are needed later on.

\begin{lemma}\label{ile}{\emph{[\cite{np}, Corollary 3.15]}}
Let $(\C,\E,\s)$ be an extriangulated category. Then the following results hold.

\item{(1)}  Let $C$ be any object in $\C$, and let
$$\xymatrix{A_{1}\ar[r]^{x_{1}}&B_{1}\ar[r]^{y_{1}}&C\ar@{-->}[r]^{\delta_{1}}&},
\xymatrix{A_{2}\ar[r]^{x_{2}}&B_{2}\ar[r]^{y_{2}}&C\ar@{-->}[r]^{\delta_{2}}&}$$
be any pair of $\E-$triangles. Then there is a commutative diagram in $\C$
$$\xymatrix{&A_{2}\ar@{=}[r]\ar[d]^{m_{2}}&A_{2}\ar[d]^{x_{2}}\\
A_{1}\ar[r]^{m_{1}}\ar@{=}[d]&M\ar[r]^{e_{1}}\ar[d]^{e_{2}}&B_{2}\ar[d]^{y_{2}}&\\
A_{1}\ar[r]^{x_{1}}&B_{1}\ar[r]^{y_{1}}&C&\\
&&}$$ which satisfies
$$\s(y_{2}^{*}\delta_{1})=[\xymatrix@C=0.41cm{A_{1}\ar[r]^{m_{1}} & M \ar[r]^{e_{1}} & B_{2}}],$$
$$\s(y_{1}^{*}\delta_{2})=[\xymatrix@C=0.41cm{A_{2}\ar[r]^{m_{2}}&M\ar[r]^{e_{2}}&B_{1}}],$$
$$m_{1*}\delta_{1}+m_{2*}\delta_{2}=0.$$
\item{(2)}  Dual of {(1)}.
\end{lemma}

The following lemma was proved in [\cite{ln}, Proposition 1.20] (see [NP, Corollary 3.16]).
\begin{lemma}\label{ilem}
Let $\xymatrix{A\ar[r]^{x}&B\ar[r]^{y}&C\ar@{-->}[r]^{\delta}&}$ be an $\E$-triangle, let $f\colon A\to D$ be any morphism, and let
$\xymatrix{D\ar[r]^{d}&E\ar[r]^{e}&C\ar@{-->}[r]^{f_{*}\delta}&}$ be any $\E$-triangle realizing $f_{*}\delta.$ Then there is a morphism $g$
which gives a morphism of $\E$-triangles
$$\xymatrix{A\ar[r]^{x}\ar[d]^f&B\ar[r]^y\ar[d]^g&C\ar@{-->}[r]^\del\ar@{=}[d]&\\
D\ar[r]^{d}&E\ar[r]^{e}&C\ar@{-->}[r]^{f_{*}\del}&}
$$
and moreover, $\xymatrix{A\ar[r]^{\binom{-f}{x}\quad}&D\oplus B\ar[r]^{\quad(d,\ g)}&E\ar@{-->}[r]^{e^{*}\delta}&}$ becomes an $\E$-triangle.
\end{lemma}

\begin{definition}{[\cite{np}, Definition 3.23]}\label{pr}
Let $\C, \E$ be as above. An object $P\in \C$ is called projective if it satisfies the following condition.
\begin{itemize}
\item For any $\E-$triangle $\xymatrix@C=0.41cm{A\ar[r]^{x} & B \ar[r]^{y}&C\ar@{-->}[r]^{\delta}&}$ and any morphism $c\in \C(P,C)$, there
    exists $b\in \C(P,B)$ satisfying $y\circ b=c$.
\end{itemize}
Injective objects are defined dually.

We denote the subcategory consisting of projective objects in $\C$ by Proj$(\C)$. Dually, the subcategory of injective objects in $\C$ is
denoted by Inj$(\C)$.
\end{definition}

\begin{lemma}{[\cite{np}, Proposition 3.24]}\label{pro}
An object $P\in \C$ is projective if and only if it satisfies $\E(P,A)=0$ for any $A\in \C$.
The dual property also holds for injective objects.
\end{lemma}

\section{Auslander-Reiten $\E-$triangles in extriangulated categories}

In this section we firstly recall the definition of a locally finite extriangulated category, which is an analogue of locally finite triangulated categories in \cite{xz}; then we will show that locally finite extriangulated categories  have Auslander-Reiten $\E-$triangles.  All extriangulated
categories in this and the next section of the paper are assumed to be Krull-Schmit, K-linear (K is a field), and dim$_{K}\mbox{Hom}_{\C}(X,Y)<\infty$, dim$_{K}\E(X,Y)<\infty$ for any
$X,Y\in\C.$

Let Ind$(\C)$ denote the set of isomorphism classes of indecomposable objects in $\C.$
For any object $A\in\C,$ denote by SuppHom$_{\C}(A,-)$ the subset of objects $B\in\C$ such that Hom$_{\C}(A,B)\neq0.$
Similarly, one can define SuppHom$_{\C}(-,A).$
If SuppHom$_{\C}(A,-)$ (SuppHom$_{\C}(-,A)$, resp.) contains only finitely many objects in Ind$(\C),$ we say that $|\mbox{SuppHom}_{\C}(A,-)|<\infty$ ($|\mbox{SuppHom}_{\C}(-,A)|<\infty,$ resp.).

\begin{definition}\label{lf}
An extriangulated category $\C$ is called locally finite if for any $A\in \mbox{Ind}(\C),$ one has $|\mbox{SuppHom}_{\C}(A,-)|<\infty$ and
$|\mbox{SuppHom}_{\C}(-,A)|<\infty.$ An extriangulated category $\C$ is called finite if $|\mbox{Ind}(\C)|<\infty.$
\end{definition}
Note that for an artin algebra $\Lambda,$ if mod$\Lambda$ is locally finite, then $|\mbox{SuppHom}_{\C}(\Lambda,-)|<\infty,$  hence $\Lambda$ is of finite representation type.

The notion of Auslander-Reiten $\E-$triangles in an extriangulated category was first introduced  in \cite{zz}, and then was fully developed in
\cite{inp}. Now we recall the definition of Auslander-Reiten $\E-$triangles in an extriangulated category.

\begin{definition}{[\cite{inp}, Definition 2.1]}\label{ar}
A non-split $\E-$extension $\delta \in\E(C,A)$ is said to be almost split if it satisfies the following conditions:
\begin{itemize}
\item (AS1) $a_{*}\delta=0$ for any non-section $a\in\C(A,A').$
\item (AS2) $c^{*}\delta=0$ for any non-retraction $c\in\C(C',C).$
\end{itemize}
\end{definition}
Note that if $a_{*}\delta=0(c^{*}\delta=0)$ for a section $a\in\C(A,A')$ (retraction $c\in\C(C',C)$), then $\delta=0.$  The realization of an
 almost split extension is called an Auslander-Reiten $\E-$triangle.

\begin{lemma}\label{elo}{\emph{[\cite{inp}, Proposition 2.5]}}
For any non-split $\E-$extension $\delta \in\E(C,A),$ the following results hold.
\begin{itemize}
 \item (1) If $\delta$ satisfies (AS1), then $A$ is indecomposable.
 \item (2) If $\delta$ satisfies (AS2), then $C$ is indecomposable.
\end{itemize}
\end{lemma}

\begin{definition}{[\cite{inp}, Definition 2.6]}\label{alm}
We say that $\C$ has Auslander-Reiten $\E-$triangles if for any indecomposable non-projective object $A\in\C,$ there exists an Auslander-Reiten
$\E-$triangle $\delta \in\E(A,B)$ for some $B\in\C$ and for any indecomposable non-injective object $B\in\C,$ there exists an
Auslander-Reiten $\E-$triangle $\delta \in\E(A,B)$ for some $A\in\C.$
\end{definition}

Note that the Auslander-Reiten sequences defined in module categories of artin algebras \cite{ars} and the Auslander-Reiten triangles defined
in triangulated categories \cite{h} are all special cases of Auslander-Reiten $\E-$triangles.

The following lemma is useful later on, which states that if we want to show a non-split $\E-$triangle is an Auslander-Reiten $\E-$triangle,
we just need to show $(AS1)$ or $(AS2)$.

\begin{lemma}\label{almo}{\emph{[\cite{inp}, Proposition 2.15]}}
Let $\delta \in\E(A,B)$ be a non-zero element such that $A$ and $B$ are indecomposable. Then the following conditions are equivalent:
\begin{itemize}
\item (1) $\delta$ is an Auslander-Reiten $\E-$triangle;
\item (2) (AS1) holds;
\item (3) (AS2) holds.
\end{itemize}
\end{lemma}

If every $\E-$triangle in $\C$ splits, then every object is projective and injective, hence the set of Auslander-Reiten $\E-$triangles is
empty. Therefore we assume that not all objects in $\C$ are projective and injective from now on.

\begin{definition}\label{def}
Let $\C$ be an extriangulated category and $X\in \mbox{Ind}(\C)$. We define a set of $\E-$triangles $S(X)$ as follows:
$$S(X)=\{s: \xymatrix{A_{s}\ar[r]^{f_{s}}&B_{s}\ar[r]^{g_{s}}&X\ar@{-->}[r]^{\delta_{s}}&}|s \mbox{\ is \ a \ non-split} \ \E-\mbox{triangle} \  \mbox{with} \ A_{s}\in \mbox{Ind}(\C)\}.$$

Dually, we can define a set of $\E-$triangles $S'(X)$ as follows:
$$S'(X)=\{s: \xymatrix{X\ar[r]^{f_{s}}&B_{s}\ar[r]^{g_{s}}&A_{s}\ar@{-->}[r]^{\delta_{s}}&}|s \mbox{\ is \ a \ non-split} \ \E-\mbox{triangle} \  \mbox{with} \ A_{s}\in \mbox{Ind}(\C)\}.$$
\end{definition}

\begin{lemma}\label{lem}
If $X\in \mbox{Ind}(\C)$ is a non-projective object, then $S(X)$ is nonempty. Dually, if $X\in \mbox{Ind}(\C)$ is a non-injective object, then $S'(X)$ is
nonempty.
\end{lemma}

\proof Since $X\in \mbox{Ind}(\C)$ is a non-projective object, there is a non-split $\E-$triangle
$\xymatrix{A\ar[r]^{f}&B\ar[r]^{g}&X\ar@{-->}[r]^{\delta}&}$ in $\C$ which ends at $X.$
Decompose $A$ into a direct sum of indecomposable objects as $A=\bigoplus\limits_{i=1}^{n} A_{i}.$
For simplicity, we can assume that $A=A_{1}\oplus A_{2}.$
Then we have the following commutative diagrams:
$$\xymatrix{A_{i}\ar[r]\ar@{=}[d]&A_{1}\oplus A_{2}\ar[r]\ar[d]^{(f_{1}\,f_{2})=f}&A_{j}\ar[d]^{a_{j}}\\
A_{i}\ar[r]&B\ar[d]^{g}\ar[r]^{h_{j}}&B_{j}\ar[d]^{b_{j}}\\
&X\ar@{=}[r]&X}$$
with $\{i, j\}=\{1, 2\}$ and $h_{1}f_{1}=a_{1}, h_{2}f_{2}=a_{2}, h_{1}f_{2}=h_{2}f_{1}=0.$
We state that at least one of the $\E-$triangles $\xymatrix{A_{j}\ar[r]^{a_{j}}&B_{j}\ar[r]^{b_{j}}&X\ar@{-->}[r]^{\delta_{j}}&}, i=1, 2$ is
non-split.
If not, using (ET3) we have the following commutative diagram:
\begin{equation}\label{t2}
\begin{array}{l}
\xymatrix{A_{1}\oplus A_{2}\ar[r]^{(f_{1}\,f_{2})=f}\ar@{=}[d]&B\ar[r]\ar[d]^{\binom{h_{1}}{h_{2}}}&X\ar[d]^{\exists\alpha}&\\
A_{1}\oplus A_{2}\ar[r]^{\left(\begin{array}{cc}a_{1}&0\\
0&a_{2}\end{array}\right)}&B_{1}\oplus B_{2}\ar[r]&X\oplus X&.}
\end{array}
\end{equation}
which shows that $\alpha^{*}0=\delta,$  contradicting to the assumption that $\delta$ is nonzero. Therefore $S(X)$ is nonempty.

Similarly, one can show that if $X\in \mbox{Ind}(\C)$ is a non-injective object, then $S'(X)$ is nonempty.\qed

\begin{definition}\label{ord}
Let $s, t\in S(X).$ We say that $s>t$ if there is an $f\in \C(A_{s},A_{t})$ such that $f_{*}s=t,$ i.e. we have the following commutative
diagram:
\begin{equation}\label{t2}
\begin{array}{l}
\xymatrix{A_{s}\ar[r]\ar[d]^{f}&B_{s}\ar[r]\ar[d]&X\ar@{=}[d]&\\
A_{t}\ar[r]&B_{t}\ar[r]&X&.}
\end{array}
\end{equation}

We say that $s\sim t$ if $f$ is an isomorphism.

Dually, let $s, t\in S'(X).$ We say that $s>t$ if there is an $f\in \C(A_{t},A_{s})$ such that $f^{*}s=t,$ i.e. we have the following commutative diagram:
\begin{equation}\label{t2}
\begin{array}{l}
\xymatrix{X\ar[r]\ar@{=}[d]&B_{t}\ar[r]\ar[d]&A_{t}\ar[d]^{f}&\\
X\ar[r]&B_{s}\ar[r]&A_{s}&.}
\end{array}
\end{equation}

We say that $s\sim t$ if $f$ is an isomorphism.

\end{definition}

\begin{lemma}\label{orde}
S(X) is a direct ordered set with the relation defined above. S'(X) is a direct ordered set with the relation defined above.
\end{lemma}

\proof Firstly, we show that if $s>t$ and $t>s,$ then $s\sim t.$

By the definition above, we have $f\in\C(A_{s},A_{t}), g\in\C(A_{t},A_{s})$ such that $f_{*}s=t$ and $g_{*}t=s.$
Let $h=gf,$ then we have $h_{*}s=s.$
Since $A_{s}$ and $A_{t}$ are indecomposable, $h$ must be an isomorphism.
If not, $h\in \mbox{rad}(\mbox{End}_{\C}(A_{s}))$ (the radical of $\mbox{End}_{\C}(A_{s})$) and $(h^{n})_{*}s=s, \forall n\geq1.$
Then $s=0$ since $\mbox{End}_{\C}(A_{s})$ is a finite dimensional local algebra, which contradicts to the fact that $s\in S(X).$

It is easy to see that if $s>t, t>u,$ we have $s>u$.

Secondly, we show that for any $s, t\in S(X),$ there is an $u\in S(X)$ such that $s>u, t>u.$

Using Lemma \ref{ile} we have the following commutative diagram:
$$\xymatrix{&A_{t}\ar@{=}[r]\ar[d]&A_{t}\ar[d]^{f_{t}}\\
A_{s}\ar[r]\ar@{=}[d]&M\ar[r]\ar[d]&B_{t}\ar[d]^{g_{t}}&\\
A_{s}\ar[r]^{f_{s}}&B_{s}\ar[r]^{g_{s}}&X&\\
&&}$$
Then by [\cite{ln}, Proposition 1.20] we have an $\E-$triangle $\xymatrix{M\ar[r]&B_{t}\oplus
B_{s}\ar[r]^{(g_{t}\,g_{s})}&X\ar@{-->}[r]^{\delta}&}$ in $\C.$
It is easy to see that $(g_{t}\,g_{s})$ is not a splitting epimorphism.
If not, there is a map $a=\binom{a_{t}}{a_{s}}\colon X\to B_{t}\oplus B_{s}$ such that
$(g_{t}\,g_{s})\cdot\binom{a_{t}}{a_{s}}=g_{t}a_{t}+g_{s}a_{s}=1_{X}.$
Since $\mbox{End}_{\C}(X)$ is local, We have that $g_{t}a_{t}$ or $g_{s}a_{s}$ is an automorphism on $X,$ which is a contradiction.

Decomposing $M$ into a direct sum of indecomposable objects as $M=\bigoplus\limits_{i=1}^{n}M_{i},$  we have the following commutative diagrams:
$$\xymatrix{\bigoplus\limits_{j\neq i}M_{j}\ar[r]\ar@{=}[d]&\bigoplus\limits_{i=1}^{n}M_{i}\ar[r]\ar[d]&M_{i}\ar[d]\\
\bigoplus\limits_{j\neq i}M_{j}\ar[r]&B_{t}\oplus B_{s}\ar[d]\ar[r]&B_{i}\ar[d]\\
&X\ar@{=}[r]&X}$$
Using the same argument as in Lemma \ref{lem}, there are $\E-$triangles $\xymatrix{M_{i}\ar[r]&B_{i}\ar[r]&X\ar@{-->}[r]^{\delta_{i}}&},
i=1,\cdots, n$\ in $\C,$ and at least one of these $\E-$triangles is non-split.
We assume that $\xymatrix{M_{1}\ar[r]&B_{1}\ar[r]&X\ar@{-->}[r]^{\delta_{1}}&}$ is non-split. Then we have the following commutative diagram:
$$\xymatrix{A_{k}\ar[r]\ar[d]^{\exists \alpha}&B_{k}\ar[r]\ar[d]&X\ar@{=}[d]\\
\bigoplus\limits_{i=1}^{n}M_{i}\ar[r]\ar[d]&B_{t}\oplus B_{s}\ar[d]\ar[r]&X\ar@{=}[d]\\
M_{1}\ar[r]&B_{1}\ar[r]&X}$$
with $k\in\{s, t\}.$ Let $u=\xymatrix{M_{1}\ar[r]&B_{1}\ar[r]&X\ar@{-->}[r]^{\delta_{1}}&},$ it is easy to see that $u$ is the desired
$\E-$triangle such that $s>u, t>u.$ Hence S(X) is a direct ordered set with the relation defined above.

Similarly, we can see that S'(X) is a direct ordered set with the relation defined above.\qed

Note that for any $X\in \mbox{Ind}(\C),$ if the minimal element exists in $S(X)$ or $S'(X),$ it is unique up to isomorphism.

\begin{lemma}\label{ilemma}
Let $\C$ be a locally finite extriangulated category. If $X\in \mbox{Ind}(\C)$ is a non-projective object, then $S(X)$ has a minimal element. Dually, if $X\in \mbox{Ind}(\C)$ is a non-injective object, then $S'(X)$ has a minimal element.
\end{lemma}

\proof Let $X\in \mbox{Ind}(\C)$ be a non-projective object.
Using Lemma \ref{lem}, there is a non-split $\E-$triangle
$$\xymatrix{A\ar[r]&\bigoplus\limits_{i=1}^{n}B_{i}\ar[r]^{(g_{1},\cdots,g_{n})}&X\ar@{-->}[r]^{\delta'}&}$$
in $\C$ with $g_{i}\in \mbox{rad}(B_{i},X).$

Since $\C$ is locally finite, there are only finite many objects $X_{i}\in \mbox{Ind}(\C), i=1,\cdots, m$ such that $\C(X_{i},X)\neq0.$
Let $f_{ij}, 1\leq j\leq d_{i}$ form a basis of the $K-$vector space $\mbox{rad}(X_{i},X).$
We set $E=(\bigoplus\limits_{i=1}^{n}B_{i})\oplus(\bigoplus\limits_{i=1}^{m}(X_{i})^{\oplus d_{i}})$ for simplicity.
Then using [\cite{ln}, Proposition, 1.20], we have an $\E-$triangle $$\xymatrix{M\ar[r]&E\ar[r]^{\alpha}&X\ar@{-->}[r]^{\delta}&}$$
where $\alpha=(g_{1},\cdots,g_{n}, f_{11},\cdots,f_{ij},\cdots,f_{md_{m}})\in \mbox{rad}(E,X).$
Hence $\delta$ is a non-split $\E-$triangle.

Decomposing $M$ into a direct sum of indecomposable objects $M=\bigoplus\limits_{i=1}^{n}M_{i}$ and using the same argument in Lemma \ref{orde},
we have a non-split $\E-$triangle $s=\xymatrix{M_{1}\ar[r]&E_{1}\ar[r]&X\ar@{-->}[r]^{\delta''}&}\in S(X)$ such that the following diagram
commutates:
$$\xymatrix{\bigoplus\limits_{i=1}^{n}M_{i}\ar[r]\ar[d]^{\eta}&E\ar[d]\ar[r]&X\ar@{=}[d]\\
M_{1}\ar[r]&E_{1}\ar[r]&X}$$

Now to show that for any non-split $\E-$triangle $\xymatrix{L\ar[r]&\bigoplus\limits_{i=1}^{m}N_{i}\ar[r]^{h}&X\ar@{-->}[r]^{\theta}&}$ in
$\C$ with $h=(h_{1},\cdots,h_{m}),$ there is a morphism $\beta \colon \bigoplus\limits_{i=1}^{m}N_{i}\to E$ such that $h=\alpha\beta.$
It is easy to see that $h\in \mbox{rad}(\bigoplus\limits_{i=1}^{m}N_{i},X)$ since $\theta$ is non-split.
Then from the definition of $f_{ij}$ and $\alpha,$ one has that such $\beta$ must exits.
Therefore, using (ET3)$\op$ we have the following commutative diagram:
$$\xymatrix{L\ar[r]\ar[d]^{\exists\gamma}&\bigoplus\limits_{i=1}^{m}N_{i}\ar[r]^{h}\ar[d]^{\beta}&X\ar@{=}[d]&\\
M\ar[r]&E\ar[r]^{\alpha}&X&.}$$
In particular, if $t=\xymatrix{L\ar[r]&\bigoplus\limits_{i=1}^{m}N_{i}\ar[r]^{h}&X\ar@{-->}[r]^{\theta}&}\in S(X),$ i.e. $t$ is a non-split
$\E-$triangle with $L\in \mbox{Ind}(\C),$ we have the following commutative diagram:
$$\xymatrix{L\ar[r]\ar[d]^{\eta\gamma}&\bigoplus\limits_{i=1}^{m}N_{i}\ar[r]^{h}\ar[d]^{\beta}&X\ar@{=}[d]&\\
M_{1}\ar[r]&E_{1}\ar[r]^{\alpha}&X&.}$$
Which shows that $s=\xymatrix{M_{1}\ar[r]&E_{1}\ar[r]&X\ar@{-->}[r]^{\delta''}&}$ is a minimal element in $S(X).$

Similarly, one can see that the dual result also holds for $S'(X).$\qed

\begin{theorem}\label{art}
Let $\C$ be a locally finite extriangulated category. If $X\in \mbox{Ind}(\C)$ is a non-projective object, there is an Auslander-Reiten $\E-$triangle
ending at $X.$
Dually, if $X\in \mbox{Ind}(\C)$ is a non-injective object, there is an Auslander-Reiten $\E-$triangle starting at $X.$ Hence $\C$ has Auslander-Reiten $\E-$triangles.
\end{theorem}

\proof Since $X\in\C$ is a non-projective object, $S(X)$ is non-empty by Lemma \ref{lem}.
Using Lemma \ref{ilemma}, there is an $\E-$triangle $s=\xymatrix{A_{s}\ar[r]^{f_{s}}&B_{s}\ar[r]^{g_{s}}&X\ar@{-->}[r]^{\delta_{s}}&}$ which is a minimal element in $S(X)$. We will show $s$ is a Auslander-Reiten $\E-$triangle ending at $X$.

Now to show that (AR2) holds for $s.$
Suppose $f\colon N\to X$ is not a retraction.
Using Lemma \ref{ilem}, we have an $\E-$triangle $\xymatrix{M\ar[r]&B_{s}\oplus N\ar[r]^{(g_{s},\ f)}&X\ar@{-->}[r]^{\theta}&},$ which is
non-split since $g_{s}, f$ are not retractions and $X$ is indecomposable.
Decomposing $M$ into a direct sum of indecomposable objects $M=\bigoplus\limits_{i=1}^{n}M_{i}$ and using the same argument as in Lemma
\ref{orde}, we have a non-split $\E-$triangle in $\C,$ and we denote it by
$u=\xymatrix{M_{1}\ar[r]&B_{1}\ar[r]&X\ar@{-->}[r]^{\delta_{1}}&}\in S(X).$
Then we have the following commutative diagram:
$$\xymatrix{A_{s}\ar[r]\ar[d]^{\exists \alpha}&B_{s}\ar[r]\ar[d]&X\ar@{=}[d]\\
\bigoplus\limits_{i=1}^{n}M_{i}\ar[r]\ar[d]&B_{s}\oplus N\ar[d]\ar[r]&X\ar@{=}[d]\\
M_{1}\ar[r]&B_{1}\ar[r]&X}$$
One can see $s>u,$ hence $s\sim u$ because $s\in S(X)$ is the minimal element in $S(X).$
In particular, we have the following commutative diagram:
$$\xymatrix{\bigoplus\limits_{i=1}^{n}M_{i}\ar[r]\ar[d]^{a}&B_{s}\oplus N\ar[d]^{(b_{1},\ b_{2})}\ar[r]^{(g_{s},\ f)}&X\ar@{=}[d]\\
A_{s}\ar[r]^{f_{s}}&B_{s}\ar[r]^{g_{s}}&X}$$
Therefore $f=g_{s}b_{2},$ and (AR2) holds for $s.$
Since $A_{s}\in \mbox{Ind}(\C),$ using Lemma \ref{almo}, we have that $s$ is an Auslander-Reiten $\E-$triangle in $\C.$

Similarly, one can see that if $X\in\C$ is a non-injective object, there is an Auslander-Reiten $\E-$triangle starting at $X.$

Hence $\C$ has Auslander-Reiten $\E-$triangles.\qed

\begin{remark}\label{exa}

If $\C$ is a locally finite triangulated category, using Theorem \ref{art}, one can see that $\C$ has Auslander-Reiten triangles (\cite{xz}, Proposition 1.3).

\end{remark}

\section{Relations for Grothendieck groups}
Given an extriangulated category $\C$ which is assumed to be Krull-Schmit, K-linear (K is a field), and satisfy that dim$_{K}\mbox{Hom}_{\C}(X,Y)<\infty$, dim$_{K}\E(X,Y)<\infty$ for any
$X,Y\in\C.$ Let $F$ be the free abelian group with a basis the set of isomorphism classes of objects in $\C$ and
$K_{0}(\C,0)$ the quotient group of $F$ by the subgroup $R,$ which is generated by elements of the form $[A]+[B]-[A\oplus B].$
It is easy to see that $K_{0}(\C,0)$ is the free abelian group generated by the isoclasses of indecomposable objects of $\C.$
The Grothendieck group $K_{0}(\C)$ of $\C$ is the factor group of $K_{0}(\C,0)$ modulo the subgroup generated by elements of the form
$[\delta]=[A]+[C]-[B]$ whenever there is an $\E-$triangle $\xymatrix{A\ar[r]^{f}&B\ar[r]^{g}&C\ar@{-->}[r]^{\delta}&}$ in $\C.$
Therefore, we have a canonical epimorphism $\psi:K_{0}(\C,0)\to K_{0}(\C).$

\begin{lemma}\label{rad}
Let $\C$ be a locally finite extriangulated category. Then for any object $A\in \C,$ there is a positive integer $n$ ($m,$ resp.) such that
$\mbox{rad}^{n}(-,A)=0$ ($\mbox{rad}^{m}(A,-)=0,$ resp.), where rad(-,-) is the radical of $\C.$
\end{lemma}

\proof The proof is the same as that in [\cite{xz}, Lemma 1.2].\qed

The following lemma is important for our investigation, which is an extriangulated category version of the result in [\cite{xz}, Lemma, 2.3].

\begin{lemma}\label{ilemm}
Let $\xymatrix{A\ar[r]^{f}&B\ar[r]^{\binom{g_{1}}{g_{2}}}&C_{1}\oplus C_{2}\ar@{-->}[r]^{\delta}&} (*)$ be an $\E-$triangle in $\C.$ If
$g_{1}$ is a retraction, then $\delta$ is isomorphic to the following $\E-$triangle:
$$\xymatrix{A\ar[r]^{\binom{0}{f_{2}}}&C_{1}\oplus B_{2}\ar[r]^{\left(\begin{array}{cc}1&0\\
0&g_{22}\end{array}\right)}&C_{1}\oplus C_{2}\ar@{-->}[r]^{\delta'}&}. (**)$$
Moreover, the $\E-$triangle $(*)$ is the direct sum of $\E-$triangles \\
$$\xymatrix{A\ar[r]^{f_{2}}&B_{2}\ar[r]^{g_{22}}&C_{2}\ar@{-->}[r]^{\delta_{1}}&} and \xymatrix{0\ar[r]&C_{1}\ar[r]&C_{1}\ar@{-->}[r]^{0}&}.$$
\end{lemma}

\proof Using (ET4)$\op$, we have the following commutative diagram:
$$\xymatrix{A\ar[r]\ar@{=}[d]&B_{2}\ar[r]\ar[d]^{h_{1}}&C_{2}\ar[d]^{\binom{0}{1}}\\
A\ar[r]^{f}&B\ar[d]^{g_{1}}\ar[r]^{\binom{g_{1}}{g_{2}}}&C_{1}\oplus C_{2}\ar[d]^{(1\ 0)}\\
&C_{1}\ar@{=}[r]&C_{1}}$$
Since $g_{1}$ is a retraction, the $\E-$triangle $\xymatrix{B_{2}\ar[r]^{h_{1}}&B\ar[r]^{g_{1}}&C_{1}\ar@{-->}[r]^{\delta'}&}$ splits.
Without losing generality, one can set $B=C_{1}\oplus B_{2}$ and $g_{1}=(1\ 0).$
Therefore  the $\E-$triangle $(*)$ is isomorphic to the following $\E-$triangle:
$$\xymatrix{A\ar[r]^{\binom{0}{f_{2}}}&C_{1}\oplus B_{2}\ar[r]^{\left(\begin{array}{cc}1&0\\
*&g_{22}\end{array}\right)}&C_{1}\oplus C_{2}\ar@{-->}[r]^{\delta''}&}. (***)$$
It is easy to see that $(***)$ is isomorphic to $(**).$

Using (ET4)$\op$, we have the following commutative diagram:
$$\xymatrix{A\ar[r]^{a}\ar@{=}[d]&B_{2}\ar[r]^{b}\ar[d]^{\alpha}&C_{2}\ar[d]^{\binom{0}{1}}\\
A\ar[r]^{\binom{0}{f_{2}}}&C_{1}\oplus B_{2}\ar[d]^{(1\ 0)}\ar[r]^{\left(\begin{array}{cc}1&0\\
0&g_{22}\end{array}\right)}&C_{1}\oplus C_{2}\ar[d]^{(1\ 0)}\\
&C_{1}\ar@{=}[r]&C_{1}}$$
Without losing generality, we can assume that $\alpha=\binom{0}{1}, a=f_{2}, b=g_{22}.$
Therefore,
$$\xymatrix{A\ar[r]^{f_{2}}&B_{2}\ar[r]^{g_{22}}&C_{2}\ar@{-->}[r]^{\delta_{1}}&} and \xymatrix{0\ar[r]&C_{1}\ar[r]&C_{1}\ar@{-->}[r]^{0}&}.$$
are $\E-$triangles and the last assertion holds.\qed

\begin{remark}\label{co}
Assume that $C\in \mbox{Ind}(\C)$ in Lemma \ref{ilemm}, then one can see that any non-split $\E-$triangle ending at $C$ is of the form
$\xymatrix{A\ar[r]&\bigoplus\limits_{i=1}^{n}B_{i}\ar[r]^{(g_{1},\cdots,g_{n})}&C\ar@{-->}[r]^{\delta'}&}$ with $g_{i}\in \mbox{rad}(B_{i},C).$
\end{remark}

\begin{theorem}\label{arr}
Let $\C$ be a locally finite extriangulated category. Then $Ker(\psi)$ is generated by elements $[\delta]$ in $K_{0}(\C,0),$ where
$\xymatrix{A\ar[r]^{f}&B\ar[r]^{g}&C\ar@{-->}[r]^{\delta}&}$ runs through all Auslander-Reiten $\E-$triangles in $\C.$
\end{theorem}

\proof It is obvious that $[A]+[C]-[B]\in \mbox{Ker}(\psi)$ if $\xymatrix{A\ar[r]^{f}&B\ar[r]^{g}&C\ar@{-->}[r]^{\delta}&}$ is an Auslander-Reiten
$\E-$triangle in $\C.$

Conversely, let $\xymatrix{A\ar[r]^{f}&B\ar[r]^{g}&C\ar@{-->}[r]^{\delta}&}$ be an arbitrary non-split $\E-$triangle in $\C.$
It suffices to prove that $[\delta]\in K_{0}(\C,0)$ can be written as a sum of $[\delta_{i}], i=1,\cdots, n,$ where each $[\delta_{i}]$
correspondences to an Auslander-Reiten $\E-$triangle.
We can assume that $C\in \mbox{Ind}(\C).$
If not, without loss of generality, we assume that $C=C_{1}\oplus C_{2}$ with $C_{1}, C_{2}\in \mbox{Ind}(\C).$
Then using (ET4)$\op$ we have the
following commutative diagram:
$$\xymatrix{A\ar[r]\ar@{=}[d]&B_{2}\ar[r]\ar[d]&C_{2}\ar[d]\\
A\ar[r]^{f}&B\ar[d]\ar[r]^{g}&C_{1}\oplus C_{2}\ar[d]\\
&C_{1}\ar@{=}[r]&C_{1}}$$
Hence $[\delta]=[\delta_{1}]+[\delta_{2}],$ where $[\delta_{1}]$ and $[\delta_{2}]$ are given by the $\E-$triangles
$\xymatrix{B_{2}\ar[r]&B\ar[r]&C_{1}\ar@{-->}[r]^{\delta_{1}}&}$ and $\xymatrix{A\ar[r]&B_{2}\ar[r]&C_{2}\ar@{-->}[r]^{\delta_{2}}&}$ with
$C_{1}, C_{2}\in \mbox{Ind}(\C)$.

Since $\delta$ does not split, $C\in \mbox{Ind}(\C)$ is not a projective object.
Using Theorem \ref{art}, there is an Auslander-Reiten $\E-$triangle
$\xymatrix{A_{1}\ar[r]^{x_{1}}&M\ar[r]^{y_{1}}&C\ar@{-->}[r]^{\delta_{1}}&}$ ending at $C.$
Since $\delta_{1}$ is an Auslander-Reiten $\E-$triangle, $y_{1}\in \mbox{rad}(M,C)$ and we have the following commutative diagram:
$$\xymatrix{A\ar[r]^{f}\ar[d]^{\alpha}&B\ar[r]^{g}\ar[d]^{\beta}&C\ar@{=}[d]&\\
A_{1}\ar[r]^{x_{1}}&M\ar[r]^{y_{1}}&C&.}$$
By Lemma \ref{ilem}, one can see that $\xymatrix{A\ar[r]^{\binom{-f}{\alpha}}&B\oplus M\ar[r]^{(\beta\
x_{1})}&M\ar@{-->}[r]^{\theta=y_{1}^{*}\delta}&}$ is an $\E-$triangle and $[\delta]=[\delta_{1}]+[\theta]$ in $K_{0}(\C,0).$
Without loss of generality, we can assume that $M=M_{1}\oplus M_{2}$ with $M_{1}, M_{2}\in \mbox{Ind}(\C).$
Assume that $B_{1}=B\oplus M.$
Then $\theta$ can be written as:
$$\xymatrix{A\ar[r]^{\psi}&B_{1}\ar[r]^{\binom{g_{1}}{g_{2}}}&M_{1}\oplus M_{2}\ar@{-->}[r]^{\theta}&}$$
with $\binom{g_{1}}{g_{2}}=(\beta\ x_{1})$ and $\psi=\binom{-f}{\alpha}.$
If $g_{1}$ is a retraction, using Lemma \ref{ilemm}, one can decompose $\theta$ into a direct sum of $\E-$triangles:
$$\xymatrix{A\ar[r]^{\psi_{1}}&B_{2}\ar[r]^{g_{22}}&M_{2}\ar@{-->}[r]^{\delta'}&} \mbox{and}
\xymatrix{0\ar[r]&M_{1}\ar[r]&M_{1}\ar@{-->}[r]^{0}&}$$
and it is easy to see that $[\theta]=[\delta'].$

If in addition, $g_{22}$ is a retraction, one can see $[\delta']=0$ and $[\theta]=[\delta']=0,$ hence $[\delta]=[\delta_{1}].$
We have proved the assertion.

If $g_{22}$ is not a retraction, then $g_{22}\in \mbox{rad}(B_{2},M_{2})$ and we have the following commutative diagram:
\begin{equation}
\begin{array}{l}
\xymatrix{A\ar[r]\ar@{=}[d]&B_{2}\ar[r]\ar[d]&M_{2}\ar[d]^{\binom{0}{1}}&\\
A\ar[r]&B_{1}\oplus B_{2}\ar[r]&M_{1}\oplus M_{2}&.}
\end{array}
\end{equation}
Let $\hat{y_{1}}=y_{1}{\binom{0}{1}},$ one can see $\hat{y_{1}}\in \mbox{rad}(M_{2},C)$ and $[\theta]=[\delta'].$
Then we can consider $\delta'$ instead of $\theta$ to continue the above process.

Now assume that $g_{1}, g_{2}$ are not retractions, hence $\binom{g_{1}}{g_{2}}$ is not a retraction.
Using Theorem \ref{art}, one has Auslander-Reiten $\E-$triangles ending at $M_{1}, M_{2}$ respectively:
$$\xymatrix{M_{i}'\ar[r]^{a_{i}}&X_{i}\ar[r]^{b_{i}}&M_{i}\ar@{-->}[r]^{\theta_{i}}&}, i=1,2.$$
Then we have the following commutative diagram:
$$\xymatrix{A\ar[r]^{\psi}\ar[d]&B_{1}\ar[r]^{\binom{g_{1}}{g_{2}}}\ar[d]&M_{1}\oplus M_{2}\ar@{=}[d]&\\
M_{1}'\oplus M_{2}'\ar[r]^{\begin{pmatrix}
a_{1} & 0 \\
0 &  a_{2}
\end{pmatrix}}&X_{1}\oplus X_{2}\ar[r]^{\begin{pmatrix}
b_{1} & 0 \\
0 &  b_{2}
\end{pmatrix}}&M_{1}\oplus M_{2}&.}$$
Let \[a=\begin{pmatrix}
a_{1} & 0 \\
0 &  a_{2}
\end{pmatrix}\] and \[b=\begin{pmatrix}
b_{1} & 0 \\
0 &  b_{2}
\end{pmatrix}.\]
Since $\theta_{i}$ are Auslander-Reiten $\E-$triangles, $b\in \mbox{rad}(X_{1}\oplus X_{2},M).$
Using Lemma \ref{ilem}, we have an $\E-$triangle
$$\xymatrix{A\ar[r]&B_{1}\oplus M_{1}\oplus M_{2}\ar[r]&X_{1}\oplus X_{2}\ar@{-->}[r]^{\theta'}&}$$
with $\theta'=b^{*}\theta=(y_{1}b)^{*}\delta,$ $y_{1}b\in \mbox{rad}^{2}(X_{1}\oplus X_{2},C)$ and $[\theta]=[\theta']+[\theta_{1}]+[\theta_{2}].$
One can continue the same process for $\theta'$ as we have done for $\theta.$
Using Lemma \ref{rad} we have $\mbox{rad}^{n}(-,C)=0$ for some $n,$ hence the process must stop at finite steps.
Therefore there are finitely many Auslander-Reiten $\E-$triangles $\delta_{i}, i=1,\cdots,n,$ such that
$[\delta]=[\delta_{1}]+\cdots+[{\delta_{n}}].$
This completes the proof.\qed

\begin{example}\label{exam}
\begin{itemize}
\item[$(1)$] Let $\Lambda$ be an artin algebra of finite representation type and mod$\Lambda$ be the category of finitely generated left $\Lambda-$modules.
Then mod$\Lambda$ is a finite extriangulated category.
Using Theorem \ref{art} and Theorem \ref{arr}, one can see the Theorem in \cite{b} holds.

\item[$(2)$] If $\C$ be a locally finite triangulated category, hence $\C$ is a locally finite extriangulated category. Using
    Theorem \ref{arr}, one can see that Theorem 2.1 in \cite{xz} holds.

\item[$(3)$] Let $\Lambda$ be a hereditary algebra and $\mathcal{P}_{0}$ be the subcategory defined by the preprojective component of the
AR quiver of $\Lambda$. It is easy to see that $\mathcal{P}_{0}$ is a locally finite exact category since $\Lambda$ is a hereditary algebra. Using Theorem \ref{arr}, we have that the relations of $K_{0}(\mathcal{P}_{0})$ are all generated by Auslander-Reiten sequences, see also \cite{mmp}. If $\Lambda$ is of infinite representation type, then the number of indecomposable objects in $\mathcal{P}_{0}$ is also infinite. In this case $\mathcal{P}_{0}$ is a locally finite category which is not finite. Note that in \cite{mmp} this example was given to show that in an exact category, the relations of the Grothendieck group are generated by Auslander-Reiten sequences doesn't imply that the exact category is finite.

\item[$(4)$] Let $\Lambda$ be a hereditary algebra and $\mathcal{P}_{0}$ and $\mathcal{I}_{0}$ be the subcategories defined by the preprojective and preinjective component of the AR quiver of $\Lambda$, let $D^{b}(\Lambda)$ be the (bounded) derived category of $\Lambda$. Then $\C=\mathcal{I}_{0}[-1]*\mathcal{P}_{0}$ is an extension closed subcategory in $D^{b}(\Lambda)$, where $\C=\mathcal{I}_{0}[-1]*\mathcal{P}_{0}=\{X\in D^{b}(\Lambda)|\exists \mbox{\ a \ triangle }\xymatrix{Y\ar[r]^{f}&X\ar[r]^{g}&Z\ar[r]&Y[1]} \ \mbox{in} \  D^{b}(\Lambda) \ \mbox{such \ that} \ Y\in \mathcal{I}_{0}[-1] \ \mbox{and} \ Z\in \mathcal{P}_{0}\}.$ Therefore, $\C$ is a extriangulated category. It is easy to see that $\C$ is locally finite since $\Lambda$ is a hereditary algebra. Therefore the relations of $K_{0}(\C)$ are all generated by Auslander-Reiten $\E-$triangles in $\C$. Note that if $\Lambda$ is of infinite representation type, then the number of indecomposable objects in $\C$ is infinite.
\end{itemize}
\end{example}

\begin{remark}\label{converse}
In general, we don't know whether the converse of Theorem \ref{arr} is true.
\end{remark}

In the rest of this subsection, we will firstly show that when restricting to triangulated categories admitting a cluster tilting subcategory $\mathcal T$, the converse of Theorem \ref{arr} is true; then we consider the relative converse of Theorem \ref{arr} in the sense of \cite{pppp} and generalize the Theorem 3.8 in \cite{pppp} to the case where $\C$ is a triangulated category with a cluster tilting subcategory $\mathcal T$. We recall the definition of cluster tilting subcategories (see also [\cite{bmrrt}, \cite{iy}, \cite{kr}, \cite{kz}]).

\begin{definition}\label{cto}
A subcategory $\mathcal T\subset\C$ which is closed under taking (finite) direct sums and direct summands is called a cluster tilting subcategory if
\begin{itemize}
\item[$(1)$] $\mathcal T$ is functorially finite in $\C;$

\item[$(2)$] $ \mathcal T=\{X\in\C|\mbox{Hom}_{\C}(\mathcal T,X[1])=0\}=\{X\in\C|\mbox{Hom}_{\C}(X,\mathcal T[1])=0\}.$

\end{itemize}
An object $T\in\C$ is called a cluster tilting object if add$(T)$ is a cluster tilting subcategory.
\end{definition}

For any $X,Y\in\C,$ one can define $[X,Y]=\mbox{dim}_{K}\mbox{Hom}_{\C}(X,Y).$ Then we have the following Lemma.

\begin{lemma}\label{ART}{\emph{[\cite{ha1}, Proposition 2.1]}}
Let $\xymatrix{X\ar[r]^{f}&Y\ar[r]^{g}&Z\ar[r]&X[1]}$ be an Auslander-Reiten triangle in $\C.$  For $A\in \mbox{Ind}(\C),$ we have that $[A,X]+[A,Z]-[A,Y]\neq0$ iff $A\cong Z$ or $A\cong Z[-1].$
\end{lemma}

\begin{theorem}\label{arr1}
Let $\mathcal T$ be a cluster tilting subcategory in $\C.$ Suppose $Ker(\psi)$ is generated by elements $[\delta]$ in $K_{0}(\C,0),$ where
$$\delta=\xymatrix{X\ar[r]^{f}&Y\ar[r]^{g}&Z\ar[r]&X[1]}$$
runs through all Auslander-Reiten triangles in $\C.$
Then $\C$ is locally finite.
In particular, if $\mathcal T=\mbox{add} T$, for a cluster tilting object $T$ in $\C,$ then $|Ind(\C)|<\infty.$
\end{theorem}

\proof Since $\mathcal T$ is a cluster tilting subcategory in $\C,$ for any nonzero indecomposable object $A\in\C,$ there is an triangle $$\xymatrix{T_{1}\ar[r]^{f}&A\ar[r]^{g}&T_{2}[1]\ar[r]&T_{1}[1]} (*)$$
with $T_{1}, T_{2} \in \mathcal T.$
Now let $T'=T_{1}\oplus T_{2}.$
For any nonzero indecomposable object $X\in\C$ with Hom$_{\C}(X,A)\neq0,$ one can see that $[X,T']+[X,T'[1]]\neq0.$
To be precise, if $[X,T']+[X,T'[1]]=0,$ we have Hom$_{\C}(X,T')=\mbox{Hom}_{\C}(X,T'[1])=0.$
Applying Hom$_{\C}(X,-)$ to $(*),$ one can see Hom$_{\C}(X,A)=0,$ contradicting to the assumption.

Considering the triangle $\theta=\xymatrix{T'\ar[r]&0\ar[r]&T'[1]\ar[r]&T'[1]},$ we have $[\theta]\in \mbox{Ker}(\psi).$
Hence $[\theta]=\sum_{i=1}^{n} a_{i}[\theta_{i}]$ with $\theta_{i}$ the AR-triangles in $\C$ such that $[X,[\theta]]=\sum_{i=1}^{n}[X,[\theta_{i}]]\neq0.$
Using Lemma \ref{ART}, one can see that for any nonzero indecomposable object $A\in\C, |\mbox{SuppHom}(-,A)|<\infty.$
Therefore, $\C$ is a locally finite triangulated category by the dual of Proposition 1.1 in \cite{xz}.

Now suppose that $T$ is a cluster tilting object in $\C.$
Consider the triangle
$$\theta=\xymatrix{T\ar[r]&0\ar[r]&T[1]\ar[r]&T[1]},$$
 we have $[\theta]\in \mbox{Ker}(\psi).$ Hence $[\theta]=\sum_{i=1}^{n} a_{i}[\theta_{i}]$ with $\theta_{i}$ the Auslander-Reiten triangles in $\C.$
 Since $T$ is a cluster tilting object in $\C,$ we have $\sum_{i=1}^{n}[A,[\theta_{i}]]=[A,[\theta]]=[A,T]+[A,T[1]]\neq0$ for any $A\in \mbox{Ind}(\C).$ Then using Lemma \ref{ART}, one can see that $|\mbox{Ind}(\C)|<\infty.$\qed

\begin{example}\label{exp}
Let $Q$ be a Dynkin quiver and mod$(KQ)$ the category of finite dimensional left $KQ-$modules.
If $T\in \mbox{mod}(KQ)$ is a tilting object, then $\mathcal T=\{F^{k}T|k\in \mathbb Z\}\subset D^{b}(KQ)$ is a cluster titling subcategory in $D^{b}(KQ)$ \cite{z}, where $F=\tau^{-1}[1]$ with $\tau$ the Auslander-Reiten translation.
Then $D^{b}(KQ)$ is a locally finite triangulated category with a culster tilting subcategory such that the relations of $K_{0}(D^{b}(KQ))$ are generated by Auslander-Reiten triangles in $D^{b}(KQ).$

\end{example}

In \cite{pppp}, the authors considered the relative Grothendieck group of a 2-CY triangulated category $\C$ admitting a cluster tilting object $T.$ They showed that the relations of the relative Grothendieck group are generated by some special Auslander-Reiten triangles iff $|\mbox{Ind}(\C)|<\infty.$
We will generalise this result to the situation where $\C$ is a triangulated category with a cluster tilting subcategory and show that the relations of the relative Grothendieck group are generated by some special Auslander-Reiten triangles iff $\C$ is locally finite.

\begin{definition}\label{gpt}
Let $\mathcal T\in\C$ be a cluster tilting subcategory. We define $K_{0}(\C,\mathcal T)$ to be the quotient of $K_{0}(\C,0)$ by the relations $[X]+[Z]-[Y]$ for all triangles
$$\xymatrix{X\ar[r]&Y\ar[r]&Z\ar[r]^{h}&X[1]}$$
with $h$ factors through an object in $\mathcal T[1].$
$K_{0}(\C,\mathcal T)$ is called the relative Grothendieck group of $\C$ with respect to $\mathcal T.$
We denote $g:K_{0}(\C,0)\to K_{0}(\C,\mathcal T)$ the canonical projection.
\end{definition}

\begin{notation}\label{lre}
Let $X\in\C$ be an indecomposable object and let
$$\xymatrix{X\ar[r]&Y\ar[r]&\tau^{-1}X\ar[r]&X[1]}$$
be an Auslander-Reiten triangle. We define
$$[l_{X}]:=[X]+[\tau^{-1}X]-[Y]\in K_{0}(\C,0).$$
\end{notation}

\begin{lemma}\label{ele}
Let $$\xymatrix{X\ar[r]&Y\ar[r]^{g}&\tau^{-1}X\ar[r]^{h}&X[1]}$$ be an Auslander-Reiten triangle. Then $X\notin \mathcal T[1]$ if and only if h factors through an object of $\mathcal  T[1].$
\end{lemma}

\proof It is easy to see that $\mathcal T[1]=\tau(\mathcal T).$
Suppose that $X\in \mathcal T[1]=\tau(\mathcal T),$ then $\tau^{-1}\X\in \mathcal T.$
If in addition $h$ factors through an object of $\mathcal  T[1],$ we have $h$=0, a contradiction.\\
Now assume that $X\notin \mathcal T[1]=\tau(\mathcal T),$ then $\tau^{-1}X\notin \mathcal T.$
Since $\mathcal T$ is a cluster tilting subcategory, we have a triangle
$$\xymatrix{T_{1}\ar[r]&T_{0}\ar[r]^{\alpha}&\tau^{-1}X\ar[r]^{\beta}&T_{1}[1]}$$
with $T_{0}, T_{1}\in\mathcal T.$
As $\tau^{-1}X\notin \mathcal T$ is indecomposable, there is a $f\in \mbox{Hom}_{\C}(T_{0},Y)$ such that $\alpha=gf.$
Hence we have $h\alpha=hgf=0,$ and then there is a $p\in \mbox{Hom}_{\C}(T_{1}[1],X[1]),$ such that $h=p\beta.$ Therefore h factors through an object of $\mathcal  T[1].$\qed

Now we have the following result.

\begin{theorem}\label{ele2}
Let $\C$ be a triangulated category which has Auslander-Reiten triangles and a cluster tilting subcategory $\mathcal T$. Then $Ker(g)$ is generated by  $\{[l_{X}]|X\in Ind(\C)\backslash \mathcal T[1]\}$ if and only if $\C$ is locally finite.
\end{theorem}

\proof Suppose that Ker$(g)$ is generated by  $\{[l_{X}]|X\in \mbox{Ind}(\C)\backslash \mathcal T[1]\}.$
For any nonzero indecomposable object $A\in\C,$ there is an triangle $$\xymatrix{T_{1}\ar[r]^{f}&A\ar[r]^{g}&T_{2}[1]\ar[r]&T_{1}[1]} $$
with $T_{1}, T_{2} \in \mathcal T.$
Now let $T'=T_{1}\oplus T_{2},$ then we have the following triangle:
$$\theta=\xymatrix{T'\ar[r]&0\ar[r]&T'[1]\ar[r]^{h}&T'[1]}.$$
It is obvious that $h$ factors through an object of $\mathcal T[1],$ we have $[\theta]\in \mbox{Ker}(g).$
Hence $[\theta]=\sum_{i=1}^{n} a_{i}[l_{X}]$ with $X\in \mbox{Ind}(\C)\backslash \mathcal T[1]$.
Using the same method in the proof of Theorem \ref{arr1}, we have that $|\mbox{suppHom}_{\C}(-,X)|<\infty.$ Therefore, $\C$ is locally finite.

Conversely, suppose that $\C$ is locally finite.
Define a subfunctor $\E_{\mathcal T}$ of Hom$_{\C}(-,-[1])$ by
$$\E_{\mathcal T}(X,Y)=\{h\in \mbox{Hom}_{\C}(X,Y[1])|h \mbox{\ factors\ through\ an\ object\ of}\ \mathcal  T[1]\}.$$
Using Proposition 3.17 in \cite{hln}, $(\C, \E_{\mathcal T})$ is an extriangulated category with all $\E-$triangles are triangles
$$\xymatrix{Y\ar[r]&Z\ar[r]&X\ar[r]^{h}&X[1]}$$
in $\C$ such that $h$ factors through an object of $\mathcal  T[1].$
Using Proposition 5.10 in \cite{inp} and Lemma \ref{ele}, one can see that $(\C, \E_{\mathcal T})$ is an extriangulated category having Auslander-Reiten $\E-$triangles since $\C$ has Auslander-Reiten triangles, and all Auslander-Reiten $\E-$triangles in $(\C, \E_{\mathcal T})$ are Auslander-Reiten triangles
$$\xymatrix{X\ar[r]&Y\ar[r]&\tau^{-1}X\ar[r]^{h}&X[1]}$$
in $\C$ such that $h$ factors through an object of $\mathcal  T[1],$ i.e. $X\notin \mathcal T[1].$
Since $\C$ is locally finite, so does $(\C, \E_{\mathcal T}).$
Using Theorem \ref{arr}, Ker$(g)$ is generated by $\{[l_{X}]|X\in \mbox{Ind}(\C)\backslash \mathcal T[1]\}.$ \qed

\begin{corollary}\label{cy1}{\emph{[\cite{pppp}, Theorem 3.8]}}
Let $\C$ be a 2-CY triangulated category with a cluster tilting object $T\in\C.$ Then $Ker(g)$ is generated by $\{[l_{X}]|X\in Ind(\C)\backslash add(T[1])\}$ if and only if $|Ind(\C)|<\infty.$
\end{corollary}

Now we have the following result.

\begin{corollary}\label{ele3}
Let $\C$ be a triangulated category which has Auslander-Reiten triangles and a cluster tilting subcategory $\mathcal T$. Then the following statements are equivalent:
\begin{itemize}
\item[$(1)$] $Ker(g)$ is generated by  $\{[l_{X}]|X\in Ind(\C)\backslash \mathcal T[1]\};$

\item[$(2)$] $Ker(\psi)$ is generated by $\{[l_{X}]|X\in Ind(\C)\};$

\item[$(3)$] $\C$ is locally finite.

\end{itemize}
\end{corollary}

\section{Classifying dense (co)resolving subcategories via Grothendieck groups}
In \cite{m}, the author used Grothendieck groups to classify dense (co)resolving subcategories in exact categories. In this section we will show
that the same results also holds for extriangulated categories, hence unifying many previous works about classifying subcategories using
Grothendieck groups.
We note that independently of our work, J. Haugland has recently generalized the result to the situation of $n-$exangulated categories in \cite{ha2}.
Recall that a subcategory $\X\subset\C$ is called dense if add$(\X)=\C.$ In this section, all category are assumed to be essentially small.

Now we recall the definition of a (co)generator in an extriangulated category.

\begin{definition}\label{ge}
Let $\mathcal{G}$ be a set of objects in $\C.$ $\mathcal{G}$ is called a generator (cogenerator, resp.) of $\C$ if for any $A\in\C,$ there is
an $\E-$triangle
$$\xymatrix{A_{1}\ar[r]&G\ar[r]&A\ar@{-->}[r]^{\delta}&} \ (\xymatrix{A\ar[r]&G\ar[r]&A_{1}\ar@{-->}[r]^{\theta}&}, resp.)$$
in $\C$ with $G\in\mathcal G.$
\end{definition}

\begin{example}\label{exap}
\begin{itemize}
\item[$(1)$] Let $\Lambda$ be an artin algebra and mod$\Lambda$ the category of finitely generated left $\Lambda-$modules.
Then $\Lambda$ is a generator of mod$\Lambda.$

\item[$(2)$] If $\C$ is a triangulated category, then $0$ is a (co)generator of $\C.$

\item[$(3)$] If $\mathcal E$ is a Frobenius category, then the subcategory $\mathcal P$ consisting of projective objects in $\mathcal E$ is
    a (co)generator of $\mathcal E.$
\end{itemize}
\end{example}

The notation of resolving subcategories were introduced by Auslander-Bridger \cite{ab}, we give the definition of $\mathcal G-$(co)resolving
subcategories of $\C$ in the setting of extriangulated categories.
\begin{definition}\label{res}
Let $\mathcal{G}$ be a set of objects and $\X$ be a subcategory in $\C.$
We say that $\X$ is a $\mathcal G-$resolving subcategory of $\C$ if the following conditions are satisfied:
\begin{itemize}
\item[$(i)$] $\X$ is closed under extensions, i.e. for any $\E-$triangle $\xymatrix{A\ar[r]&B\ar[r]&C\ar@{-->}[r]^{\delta}&}$ with
    $A,C\in\X,$ we have $B\in\X;$

\item[$(ii)$] $\X$ is closed under CoCone of deflations;

\item[$(iii)$] $\mathcal G\subset\X.$
\end{itemize}

\end{definition}
Dually, one can define $\mathcal G-$coresolving subcategories of $\C.$
\begin{definition}\label{cres}
Let $\mathcal{G}$ be a set of objects and $\X$ be a subcategory in $\C.$
We say that $\X$ is a $\mathcal G-$coresolving subcategory of $\C$ if the following conditions are satisfied:
\begin{itemize}
\item[$(i)$] $\X$ is closed under extensions, i.e. for any $\E-$triangle $\xymatrix{A\ar[r]&B\ar[r]&C\ar@{-->}[r]^{\delta}&}$ with
    $A,C\in\X,$ we have $B\in\X;$

\item[$(ii)$] $\X$ is closed under Cone of inflations;

\item[$(iii)$] $\mathcal G\subset\X.$
\end{itemize}

\end{definition}

The following two Lemmas are extriangulated category verision of Lemma 2.5 and Lemma 2.9 in \cite{m}.

\begin{lemma}\label{lm1}
Let $\X$ be a dense subcategory of $\C.$ Then $\X$ is closed under Cone of inflations iff it is closed under CoCone of deflations.
\end{lemma}

\proof The proof of the lemma is the same as that of [\cite{m}, Lemma 2.5].\qed

\begin{lemma}\label{lm2}
Let $\mathcal G$ be a generator of $\C$ and $\X$ be a dense $\mathcal G-$resolving subcategory of $\C.$ Then for any $A\in\C,$ we have $A\in\X$
iff $[A]\in \langle[X]|X\in\X\rangle,$ where $\langle[X]|X\in\X\rangle$ is the subgroup of $K_{0}(\C)$ generated by elements [X] with
$X\in\X.$
\end{lemma}

\proof It is easy to see that if $A\in\X,$ $[A]\in \langle[X]|X\in\X\rangle.$

Now to show the converse.
One can define an equivalence relation $\thicksim$ on $\C$ as follows: $C_{1}\thicksim C_{2}$ if there are $X_{1},X_{2}\in\X$ such
that $C_{1}\oplus X_{1}\cong C_{2}\oplus X_{2}.$
Let $\langle\C\rangle_{\X}:=\C/\thicksim$ and denote by $\langle C\rangle$ the class of $C.$
Then it is easy to see that $\langle\C\rangle_{\X}$ is an abelian group with $\langle A\rangle+\langle B\rangle=\langle A\oplus B\rangle.$
Let $\xymatrix{A\ar[r]&B\ar[r]&C\ar@{-->}[r]^{\delta}&}$ be an $\E-$triangle. Then there exist $A_{1},C_{1}\in\C$ such that
$A\oplus A_{1},C\oplus C_{1}\in\X$ since $\X$ is a dense subcategory in $\C.$
Hence we have the following $E-$triangle:
$$\xymatrix{A\oplus A_{1}\ar[r]&B\oplus A_{1}\oplus C_{1}\ar[r]&C\oplus C_{1}\ar@{-->}[r]^{\delta'}&}.$$
One can see that $B\oplus A_{1}\oplus C_{1}\in\X$ since $\X$ is closed under extensions.
Then we have $\langle B\rangle-\langle A\rangle-\langle C\rangle=\langle B\oplus A_{1}\oplus C_{1}\rangle-\langle A\oplus A_{1}\rangle-\langle C\oplus C_{1}\rangle=\langle 0\rangle.$
Hence we have the following group homomorphism
$$\psi: K_{0}(\C)\rightarrow \langle\C\rangle_{\X}.  \ [A]\longmapsto \langle A\rangle.$$
Note that $\langle[X]|X\in\X\rangle\subset \mbox{Ker}(\psi).$
Any element in $K_{0}(\C)$ can be written as $[A]-[B]$ with $A,B\in\C.$
Since $\mathcal G$ is a generator of $\C,$ then for $B$ we have the following $\E-$triangle
$\xymatrix{B_{1}\ar[r]&G\ar[r]&B\ar@{-->}[r]^{\delta}&}$ with $G\in\mathcal G.$
Hence $[B]=[G]-[B_{1}]$ in $K_{0}(\C),$ and $[A]-[B]=[A\oplus B_{1}]-[G],$ i.e. every element in $K_{0}(\C)$ can be written as $[A]-[G]$ with
$A\in\C$ and $G\in\mathcal G.$

Now suppose that $[A]\in \langle[X]|X\in\X\rangle.$ Then $[A]=[A_{1}]-[G]\in\langle[X]|X\in\X\rangle\subset \mbox{Ker}(\psi)$ with $G\in\mathcal G.$
Since $\mathcal G\subset \X,$ we have $[A_{1}]\in \mbox{Ker}(\psi),$ i.e. $\langle A_{1}\rangle=\langle 0\rangle.$
Hence there are $X_{1},X_{2}\in\X$ such that $A_{1}\oplus X_{1}=X_{2}.$
Consider the following split $\E-$triangle $\xymatrix{A_{1}\ar[r]&A_{1}\oplus X_{1}\ar[r]&X_{1}\ar@{-->}[r]^{\delta'}&},$ we have $A_{1}\in\X$
since $\X$ is a $\mathcal G-$resolving subcategory.
Then we have $[A]=[A_{1}]-[G]$ with $A_{1}, G\in\X,$ i.e. $[A\oplus G]=[A_{1}]$ with $A_{1}\in\X.$
Therefore, $\langle A\rangle=\langle 0\rangle.$
Using the same method as above, one can see that $A\in\X.$\qed

\begin{theorem}\label{thm}
Let $\C$ be an essentially small extriangulated category with a (co)generator $\mathcal G.$ Then we have a one-to-one correspondence between
the following sets:
\begin{itemize}
\item[$(1)$] {dense $\mathcal G-$resolving subcategories of $\C$};

\item[$(2)$] {dense $\mathcal G-$coresolving subcategories of $\C$};

\item[$(3)$] {subgroups of $K_{0}(\C)$ containing the image of $\mathcal G$}.
\end{itemize}
\end{theorem}

\proof Using Lemma \ref{lm1}, one can see that there is a one-to-one correspondence between (1) and (2).

Now to show that there is a one-to-one correspondence between (1) and (3).

For a dense $\mathcal G-$resolving subcategory $\X,$ define
$$f(\X):=\langle[X]|X\in\X\rangle$$
For a subgroup $H$ of $K_{0}(\C)$ containing the image $\mathcal G,$ define
$$g(H):=\{A\in\C|[A]\in H\}$$

Now to show that $f$ and $g$ give mutually inverse between (1) and (3).

It is obvious that $f(\X)$ is a subgroup of $K_{0}(\C)$ containing the image of $\mathcal G.$

For any subgroup $H$ of $K_{0}(\C)$ containing the image of $\mathcal G,$ to show that $g(H)$ is a dense $\mathcal G-$resolving subcategory of
$\C.$
For any $A\in\C,$ we have an $\E-$triangle $\xymatrix{A_{1}\ar[r]&G\ar[r]&A\ar@{-->}[r]^{\delta}&}$ with $G\in\mathcal G.$
Hence $[A\oplus A_{1}]=[A]+[A_{1}]=[G]\in H$ and $A\oplus A_{1}\in g(H),$ i.e. $g(H)$ is dense in $\C.$
For any $\E-$triangle $\xymatrix{A\ar[r]&B\ar[r]&C\ar@{-->}[r]^{\delta}&}$ in $\C,$ since $[A]+[C]-[B]=[0],$ we have $A\in g(H)$ if $B, C\in
g(H).$
Then one can see $g(H)$ is a dense $\mathcal G-$resolving subcategory in $\C.$

Let $H$ be a subgroup of $K_{0}(\C)$ containing $\mathcal G.$
It is easy to see that $fg(H)\subset H.$
For any $[A]-[G]\in H$ with $G\in\mathcal G,$ since $[A]=[A]-[G]+[G]\in H,$ we have $A\in g(H).$
Hence $[A]-[G]\in fg(H)$ and therefore $fg(H)=H.$

Let $\X$ be a dense $\mathcal G-$resolving subcategory of $\C.$
It is easy to see that $\X\subset gf(\X).$
For any $A\in gf(\X),$ using Lemma \ref{lm2}, we have $A\in\X$ since $[A]\subset f(\X).$
Hence $gf(\X)=\X.$ This completes the proof.\qed

\begin{corollary}\label{exp}
\begin{itemize}
\item[$(1)$]\emph[\cite{t}, Theorem 2.1] Let $\C$ be an essentially small triangulated category. Then
there is a one-to-one correspondence between the dense triangulated subcategories of $\C$ and the subgroups H of $K_{0}(\C).$

\item[$(2)$]\emph[\cite{m}, Theorem 2.7] If $\C$ be an essentially small exact category with either a generator or a cogenerator $\mathcal G$.
Then there is a one-to-one correspondence between dense $\mathcal G$-(co)resolving subcategories of $\C$ and the subgroups H of $K_{0}(\C)$ containing the image of $\mathcal G.$

\end{itemize}
\end{corollary}

\section*{Acknowlegement}

The authors would like to thank Aslak Buan, Panyue Zhou for useful comments and remarks.

Bin Zhu\\
Department of Mathematical Sciences, Tsinghua University,
100084 Beijing, P. R. China.\\
E-mail: \verb"zhu-b@mail.tsinghua.edu.cn"\\[0.3cm]
Xiao Zhuang\\
Department of Mathematical Sciences, Tsinghua University,
100084 Beijing, P. R. China.\\
E-mail: \verb"zhuangx16@mails.tsinghua.edu.cn"


\begin{thebibliography}{99}

\bibitem[A]{a} M. Auslander. Relations for Grothendieck Groups of Artin Algebras. Proceedings of the American Mathematical Society,
    91(3):336-340 1984.

\bibitem[AB]{ab} M. Auslander, M. Bridger. Stable module theory. Memoirs of the American Mathematical Society, Vol.94 (American Mathematical Society, Providence, RI, 1969).


\bibitem[AR1]{ar1} M. Auslander, I. Reiten. Representation theory of Artin algebras III, Comm. Algebra. 3:239-294, 1975.

\bibitem[AR2]{ar2} M. Auslander, I. Reiten. Grothendieck groups of algebras and orders. Journal of Pure and Applied Algebra. 39:1-51,
    1986.

\bibitem[ARS]{ars} M. Auslander, I. Reiten, S. O. Smal{\o}. Representation theory of Artin algebras. Cambridge Studies in Advanced
    Mathematics, 36. Cambridge University Press, Cambridge. 1997.


\bibitem[ASm]{asm} M. Auslander, S. O. Smal{\o}. Almost split sequences in subcategories. Journal of Algebra. 69(2):426-454, 1981.




\bibitem[B]{b} M. C. R. Butler. Grothendieck groups and almost split sequences. Integral Representations and Applications. Lecture Notes in Math.,882 Springer, Berlin-New York. 357-368, 1981.


\bibitem[Be]{be} A. Beligiannis. Auslander-Reiten Triangles, Ziegler Spectra and Gorenstein Rings. K-Theory. 32(1):1-82, 2004.


\bibitem[BMRRT]{bmrrt} A. Buan, R. Marsh, M. Reineke, I. Reiten, G. Todorov. Tilting theory and cluster combinatorics. Adv. Math. 204(2):572-618, 2006.

\bibitem[CF]{cf} B. Keller, I. Reiten. Cluster-tilted algebras are Gorenstein and stably Calabi-Yau. Adv. Math. 211(1):123-151, 2007.











\bibitem[H]{h}  D. Happel. Triangulated Categories in the Representation of Finite Dimensional Algebras. London Mathematical Society Lecture
    Note Series, 119. Cambridge University Press, 1988.


\bibitem[Ha1]{ha1} J. Haugland. Auslander-Reiten triangles and Grothendieck groups of triangulated categories. arXiv preprint arXiv:1904.02506, 2019.

\bibitem[Ha2]{ha2} J. Haugland. The Grothendieck group of an n-exangulated category. arXiv preprint arXiv:1912.04328, 2019.



\bibitem[HLN]{hln} H. Herschend, Y. Liu, H. Nakaoka. n-exangulated categories. arXiv preprint arXiv:1709.06689, 2017.

\bibitem[I]{i} O. Iyama. $\tau$-categories I: Ladders. Algebras and Representation Theory. 8(3):297-321, 2005.

\bibitem[IY]{iy} O. Iyama, Y. Yoshino. Mutations in triangulated categories and rigid Cohen-Macaulay modules. Invent. Math. 172(1):117-168, 2008.

\bibitem[INP]{inp} O. Iyama, H. Nakaoka, Y. Palu. Auslander-Reiten theory in extriangulated categories. arXiv preprint arXiv:1805.03776,
    2018.

\bibitem[J]{j} P. J{\o}rgensen. Auslander-Reiten triangles in subcategories. Journal of K-theory: K-theory and its Applications to Algebra, Geometry, and Topology. 3(3):583-601, 2009.


\bibitem[K]{k} M, Kleiner. Approximations and almost split sequences in homologically finite subcategories. Journal of Algebra. 198(1):135-163 1997.

\bibitem[KR]{kr} B. Keller, I. Reiten. Cluster-tilted algebras are Gorenstein and stably Calabi-Yau. Advances in Mathematics. 211(1):123-151, 2007.

\bibitem[KZ]{kz} S. Koenig, B. Zhu. From triangulated categories to abelian categories: cluster tilting in a general framework. Math. Z. 258(1):143-160, 2008.

\bibitem[Liu]{liu} S. Liu. Auslander-Reiten theory in a Krull-Schmidt category. The S$\tilde{a}$o Paulo of Mathematical Sciences. 4(3):425-472 2010.

\bibitem[LN]{ln} Y. Liu, H. Nakaoka. Hearts of twin Cotorsion pairs on extriangulated categories. Journal of Algebra. 528:96-149, 2019.

\bibitem[M]{m} H . Matsui. Classifying dense subcategories of exact categories via Grothendieck groups. Algebras and Representation Theory.
    21(7):1-13, 2016.

\bibitem[MMP]{mmp} E. N. Marcos, H. Merklen, M. I. Platzeck. The Grothendieck group of the category of modules of finite projective dimension over certain weakly triangular algebras. Communications in Algebra. 28(3):1387-1404, 2000.


\bibitem[NP]{np}  H. Nakaoka, Y. Palu. Extriangulated categories, Hovey twin cotorsion pairs and model structures. Cah. Topol. G$\acute{e}$eom. Diff$\acute{e}$r. Cat$\acute{e}$g. LX, 60(2),117-193, 2019.


\bibitem[PPPP]{pppp} A. Padrol, Y. Palu, V. Pilaud, P. G. Plamondon. Associahedra for finite type cluster algebras and minimal relations between $\mathbf {g} $-vectors. arXiv preprint arXiv:1906.06861, 2019.






\bibitem[R]{r} C. M. Ringel. Tame algebras and integral quadratic forms. Lecture Notes in Mathematics, 1099. Springer-Verlag, Berlin, 1984.


\bibitem[Sh]{sh} A. Shah. Auslander-Reiten theory in quasi-abelian and Krull-Schmidt categories. Journal of Pure and Applied Algebra. 224(1):98-124, 2020.

\bibitem[T]{t} R. W. Thomason. The classification of triangulated subcategories. Compositio Mathematica. 105(1):1-27, 1997.

\bibitem[XZ]{xz} J. Xiao, B. Zhu. Relations for the Grothendieck groups of triangulated categories. Journal of Algebra. 257(1):37-50, 2002.




\bibitem[ZZ]{zz} P. Zhou, B. Zhu. Triangulated quotient categories revisited. Journal of Algebra. 502:196-232 2018.



\bibitem[Z]{z} B. Zhu. Cluster-Tilted Algebras and Their Intermediate Coverings. Communications in Algebra, 39(7):2437-2448 2011.
































\end{thebibliography}
\end{document}